\documentclass[reqno,11pt,a4paper]{amsart}
\usepackage{graphics}
\usepackage{amssymb,latexsym,amsfonts,amsmath}
\usepackage{mathrsfs}
\newtheorem{lemma}{\bf Lemma}[section]

\newtheorem{proposition}{\bf Proposition}[section]

\newtheorem{theorem}{\bf Theorem}[section]

\numberwithin{equation}{section}
\theoremstyle{plain}
\theoremstyle{definition}

\newtheorem{remark}{Remark}[section]
\newtheorem{example}{Example}[section]

\newcommand{\eps}{\varepsilon}
\newcommand{\Real}{\mathbb{R}}
\renewcommand{\P}{\mathsf{P}}
\newcommand{\E}{\mathsf{E}}
\DeclareMathOperator*{\trace}{trace}
\newcommand{\lef}{\langle\hskip -1.8pt \langle}
\newcommand{\rig}{\rangle\hskip -1.8pt \rangle}
\renewcommand{\b}[1]{\mbox{\boldmath $#1$}}

\begin{document}

\title[]{The Freidlin-Wentzell LDP with rapidly growing coefficients}
\author{P. Chigansky}
\address{Department of Mathematics, The Weizmann Institute of Science,
Rehovot 76100, Israel} \email{pavel.chigansky@weizmann.ac.il}
\thanks{Research supported by a grant from the Israel Science Foundation}

\author{R. Liptser}
\address{Department of Electrical Engineering Systems,
Tel Aviv University, 69978 Tel Aviv, Israel}
\email{liptser@eng.tau.ac.il}

%\subjclass{}%
%\commby{}%
% ----------------------------------------------------------------

\begin{abstract}
The Large Deviations Principle (LDP) is verified for a homogeneous
diffusion process with respect to a Brownian motion $B_t$,
$$
X^\eps_t=x_0+\int_0^tb(X^\eps_s)ds+
\eps\int_0^t\sigma(X^\eps_s)dB_s,
$$
where  $b(x)$ and $\sigma(x)$ are are locally Lipschitz functions
with super linear growth. We assume that the drift is directed towards
the origin and the growth rates of the drift and diffusion terms are properly
balanced. Nonsingularity of $a=\sigma\sigma^*(x)$ is not required.
\end{abstract}
\maketitle
% ----------------------------------------------------------------
\section{\bf Introduction}
\label{sec-1}

In this paper we extend the set of conditions, under
which Freidlin-Wentzell's Large Deviation Principle (LDP) for a
homogeneous diffusion process remains valid. We consider a family
$\{(X^\eps_t)_{t\ge 0}\}_{\eps\to 0}$ of diffusions,
where $X^\eps_t\in\Real^d$, $d\ge 1$ is defined by the
It\^o equation
\begin{equation}\label{1.1}
X^\eps_t=x_0+\int_0^tb(X^\eps_s)ds+
\eps\int_0^t\sigma(X^\eps_s)dB_s,
\end{equation}
relative to a standard Brownian motion $B_t$, where $b(x)$ and
$\sigma(x)$ are vector and matrix valued continuous functions of
dimensions $d$ and $d\times d$ respectively, guaranteeing
existence of the unique weak solution.

The classical Freidlin-Wentzell setting \cite{FW} (see e.g.
Dembo and Zeitouni, \cite{DZ}) covers the model  \eqref{1.1} with
bounded $b(x)$ and $\sigma(x)$ and uniformly positive definite
diffusion matrix $a(x)=\sigma\sigma^*(x)$. Various LDP versions
can be found in Dupuis and Ellis \cite{11}, Feng \cite{14}, Feng
and Kurtz \cite{15}, Friedman \cite{19}, Liptser and Pukhalskii,
\cite{LP}, Mikami \cite{29}, Narita \cite{31}, Stroock \cite{40},
Ren and Zhang \cite{RZ}. In the recent paper
\cite{PP3}, Puhalskii extends LDP to \eqref{1.1} with continuous
and unbounded coefficients and singular $a(x)$, assuming $b(x)$ and
$a(x)$ are Lipschitz continuous functions (concerning singular
$\sigma(x)$ see also Liptser et al,
\cite{LSV}). Being  Lipschitz continuous, the entries of
$b,\sigma$ grow not faster than linearly and, thereby, {\it
automatically} guarantee one of the necessary conditions for LDP
($\|\cdot\|$ denotes the Euclidean norm in $\Real^d$)
\begin{equation}\label{CCCT}
\lim_{C\to\infty}\varlimsup_{\eps\to
0}\eps^2\log\P \Big(\sup_{t\le T}
\|X^\eps_t\|>C\Big)=-\infty, \ \forall \ T>0.
\end{equation}
Relinquishing the linear growth condition for $b,\sigma$ would
require additional assumptions providing \eqref{CCCT}.

This paper is inspired by Puhalskii's remark in
\cite{PP3}:
\begin{quote}\em
If the drift is directed towards the origin, then no restrictions
are needed on the growth rate of the drift coefficient.
\end{quote}
In particular, in this case the LDP holds, regardless of the growth rate of $b(x)$,
for a constant diffusion matrix (not necessarily nonsingular).

In this paper, we show that in fact LDP remains valid for
\eqref{1.1} with non-constant diffusion term, if its growth rate
is properly balanced  relatively to the drift (see \eqref{tri} of Theorem
\ref{theo-2.1} below). Our result is formulated in terms of Khasminskii-Veretennikov's
condition \eqref{dva} (see \cite{Kh} and
\cite{PVI}, \cite{PVII})

The rest of the paper is organized as follows. In Sections
\ref{sec-2} and \ref{sec-3}, the main result, notations and
preliminary facts on the LDP are given. Sections \ref{sec-4} -
\ref{sec-6} contain the proof of the main result. Auxiliary technical details
are gathered in Appendices  \ref{sec-A} - \ref{sec-B}.

\section{\bf Notations and the main result}
\label{sec-2}

The following notations and conventions are used through the
paper.
\begin{list}{-}{}
\item $^*$ denotes the transposition symbol

\item all vectors are columns (unless explicitly stated otherwise)

\item $|x|$ and $\|x\|$ denote the $\ell_1$ and $\ell_2$
(Euclidian) norms of $x\in \Real^d$

\item $\lef x,y\rig$ denotes the scalar product of $x,y\in
\Real^d$

\item $\|x\|^2_\Gamma=\lef x,\Gamma x\rig$ with an nonnegative
definite matrix $\Gamma$

\item $a(x)=\sigma(x)\sigma^*(x)$

\item $a^\oplus(x)$ denotes the Moore-Penrose pseudoinverse matrix
of $a(x)$ (see \cite{Al})

\item $\nabla V(x)$ is the  gradient (row) vector of $V(x)$:
$$
\nabla
V(x):=\Big(\frac{\partial V(x)}{\partial
x_1},\ldots,\frac{\partial V(x)} {\partial x_d}\Big)
$$

\item $\langle M,N\rangle_t$ is the joint quadratic variation process of
continuous martingales $M_t$ and $N_t$; for brevity $\langle M,M\rangle_t=
\langle M\rangle_t$

\item a.s. abbreviates  ``almost surely''; when the corresponding measure is not specified
the Lebesgue measure on $\Real_+$ is understood

\item $\varrho$ is the locally uniform metric on $\mathbb{C}_{[0,\infty)}
(\Real^d)$

\item $\mathbf{I}$ denotes $d\times d$ identity matrix

\item the convention $0/0=0$ is kept throughout

\item $X^\eps=(X^\eps_t)_{t\ge 0}$

\item $\inf\{\varnothing\}=\infty.$
\end{list}

\medskip
\noindent
We study the LDP for the family $\{X^\eps\}_{\eps\to 0}$
in the metric space $(\mathbb{C}_{[0,\infty)}(\Real^d),\varrho)$
with
$\varrho(x,y)=\sum_{k=1}^\infty 2^{-k}\big(1\vee \sup_{t\le k}\|x_t-y_t\|\big)$, $x,y\in \mathbb{C}_{[0,\infty)}(\Real^d)$.
Recall that  $\{X^\eps\}_{\eps\to 0}$ satisfies the
LDP with the good rate function $J(u): \mathbb{C}_{[0,\infty)}(\Real^d)\mapsto [0,\infty]$ and the
rate $\eps^2$, if the level sets of $J(u)$ are compacts and for any
closed set $F$ and open set $G$ in
$\mathbb{C}_{[0,\infty)}(\Real^d)$,

\begin{align*}
&\varlimsup_{\eps\to 0}\eps^2\log\P\big(X^\eps\in F\big)\le -\inf_{u\in F}J(u), \\
&\varliminf_{\eps\to 0}\eps^2\log\P\big(X^\eps\in G\big)\ge
-\inf_{u\in G}J(u).
\end{align*}

Our main result is
\begin{theorem}\label{theo-2.1}
Assume{\rm :}
\begin{enumerate}
\renewcommand{\theenumi}{H-\arabic{enumi}}

\item \label{raz} entries of $b(x)$ and $\sigma(x)$ are locally
Lipschitz continuous functions,

\item \label{dva} $\lim_{\|x\|\to\infty}\dfrac{\lef
x,b(x)\rig}{\|x\|}=-\infty$,

\item \label{tri} for some positive constants $K$ and $L$,
$\dfrac{\lef x,a(x)x\rig}{\|x\| \ |\lef x,b(x)\rig|}\le K$, $\forall
\ \|x\|>L$.
\end{enumerate}

\noindent Then $\{X^\eps_t\}_{\eps\to
0}$ obeys the LDP in the metric space
$(\mathbb{C}_{[0,\infty)}(\Real^d),\varrho)$ with the rate
$\eps^2$ and the rate function
$$
J(u)=
  \begin{cases}
    \frac{1}{2}\int_0^\infty \|\dot{u}_t-b(u_t)\|^2_{a^\oplus(u_t)}dt, &
u\in\varGamma
  \\
    \infty , & u\not\in\varGamma,
  \end{cases}
$$
where
$$
\varGamma=\Big\{u\in\mathbb{C}_{[0,\infty)}:
\begin{array}{c}
u_0=x_0, \ du_t\ll dt, \ \int_0^\infty\|\dot{u}_t\|^2dt<\infty \\
a(u_t)a^\oplus(u_t)[\dot{u}_t-b(u_t)]=[\dot{u}_t-b(u_t)] \ \text{a.s.}
\end{array}
\Big\}.
$$
\end{theorem}
\begin{remark}
In the scalar case (recall 0/0=0)
$$
J(u)=
  \begin{cases}
    \dfrac{1}{2}\int_0^\infty \dfrac{(\dot{u}_t-b(u_t))^2}{\sigma^2(u_t)}dt, &
    du_t=\dot{u}_tdt, \ u_0=x_0, \ \int_0^\infty\dot{u}^2_tdt<\infty
  \\
    \infty , & \text{otherwise}.
  \end{cases}
$$
\end{remark}

\begin{example}
A typical example within the scope of Theorem \ref{theo-2.1} is
$$
X^\eps_t=x_0-\int_0^t(X^\eps_s)^3ds+\eps\int_0^t
|X^\eps_s|^{3/2}dB_s.
$$
\end{example}

\section{\bf Preliminaries}
\label{sec-3}

We follow the framework, set up by  A.Puhalskii (see \cite{P1},
\cite{puh2}):
$$
\left.
\begin{array}{ll}
    \text{Exponential tightness} &
    \\
    \text{Local LDP}  &
  \end{array}
  \right\}
  \Longleftrightarrow \text{LDP}
$$
The exponential tightness in the metric space
$(\mathbb{C}_{[0,\infty)},\varrho)$ is convenient to verify in
terms of, so called, {\em $\mathbb{C}$-exponential tightness} conditions
introduced by A.Puhalskii (see e.g. \cite{LP}), which are based on D. Aldous's
``stopping time and tightness'' concept (see \cite{Ald1}, \cite{Ald2}).
To this end, let us assume that the
diffusion processes are defined on a stochastic basis
$(\varOmega,\mathcal{F},\mathbf{F}^\eps=(\mathscr{F}^\eps_t)_{t\ge
0}, \P)$, satisfying the usual conditions, where the
filtration $\mathbf{F}^\eps$ may depend on $\eps$.

Recall (see \cite{LP}) that the family of diffusion processes
is
$\mathbb{C}$-exponentially tight if for any $T>0$, $\eta>0$ and any
$\mathbf{F}^\eps$-stopping time $\theta$,
\begin{eqnarray}\label{3.1abs}
&& \lim_{C\to\infty}\varlimsup_{\eps\to
0}\eps^2\log\P \Big(\sup_{t\le T}
\|X^\eps_t\|>C\Big)=-\infty,
\\
&&\lim_{\triangle\to 0}\varlimsup_{\eps\to
0}\eps^2\log\sup_{\theta\le T} \P\Big(\sup_{t\le
\triangle}\|X^\eps_{\theta+t}-X^\eps_\theta\|>\eta
\Big)=-\infty. \label{3.1abd}
\end{eqnarray}
The family of diffusion processes obeys {\em the local LDP} in
$(\mathbb{C}_{[0,\infty)}(\Real^d),\varrho)$ if for any $T>0$ there exists a local
rate function $J_T(u)$ such that
\begin{eqnarray}\label{3.2abs}
&& \varlimsup_{\delta\to 0}\varlimsup_{\eps\to
0}\eps^2\log\P \Big(\sup_{t\le T}
\|X^\eps_t-u_t\|\le\delta\Big)\le-J_T(u)
\\
&& \varliminf_{\delta\to 0}\varliminf_{\eps\to
0}\eps^2\log\P \Big(\sup_{t\le T}
\|X^\eps_t-u_t\|\le\delta\Big)\ge-J_T(u).
\label{3.2abd}
\end{eqnarray}

Under the conditions \eqref{3.1abs}-\eqref{3.2abd}, the family of diffusion processes
obeys the LDP with the rate $\eps^2$ and the
good rate function
$$
J(u)=\sup_T J_T(u), \ u\in \mathbb{C}_{[0,\infty)}(\Real^d),
$$
where
$$
J_T(u)=
  \begin{cases}
    \frac{1}{2}\int_0^T \|\dot{u}_t-b(u_t)\|^2_{a^\oplus(u_t)}dt, &
u\in\varGamma_T
  \\
    \infty , & u\not\in\varGamma_T,
  \end{cases}
$$
with
$$
\varGamma_T=\Big\{u\in\mathbb{C}_{[0,T]}:
\begin{array}{c}
u_0=x_0, \ du_t\ll dt, \ \int_0^T\|\dot{u}_t\|^2dt<\infty \\
a(u_t)a^\oplus(u_t)[\dot{u}_t-b(u_t)]=[\dot{u}_t-b(u_t)] \ \text{a.s.}
\end{array}
\Big\}.
$$
Thus
the proof of Theorem \ref{theo-2.1} is reduced to establishing \eqref{3.1abs} -
\eqref{3.2abd}.

\section{\bf The proof of $\mathbb{C}$-exponential tightness}
\label{sec-4}
\subsection{Auxiliary lemma}

 Let $\mathfrak{D}$ be a nonlinear operator acting on
continuously differentiable functions $V(x):\Real^d\to\Real$ as
follows:
\begin{equation*}
\mathfrak{D}V(x)=\lef \nabla V(x), b(x)\rig+\frac{1}{2} \lef
\nabla V(x), a(x)\nabla V(x)\rig.
\end{equation*}
\begin{lemma}\label{lem-2.1}
Assume there exists twice continuously differentiable nonnegative
function $V(x)$ such that
\begin{enumerate}
\renewcommand{\theenumi}{a-\arabic{enumi}}
\item \label{a-1} $\lim_{C\to\infty}\inf_{\|x\|\ge C}V(x)=\infty$
\item \label{a-2} for some $L>0$, $\mathfrak{D}V(x)\le 0, \
\forall \ \|x\|>L$.
\end{enumerate}
Then \eqref{3.1abs} holds.
\end{lemma}
\begin{proof}
Notice that \eqref{3.1abs} is equivalent to
\begin{equation}\label{TC}
\lim_{C\to\infty}\varlimsup_{\eps\to
0}\eps^2\log\P\big( \varTheta_C\le T\big) =-\infty,
\end{equation}
where
\begin{equation}\label{ttt}
\varTheta_C=\inf\{t:\|X^\eps_t\|\ge C\}, \quad C>0
\end{equation}
are stopping times relative to $\mathbf{F}^\eps$.

We use \eqref{PA-1} of Proposition \ref{pro-A.1} to estimate
$
\log\P(\varTheta_C\le T).
$
An appropriate martingale $M^\eps_t$ is constructed with the
help of
function $V(x)$. Let $V_{ij}(x)=\frac{\partial^2
V(x)}{\partial x_i\partial x_j}$ and define
$$
\Psi(x)=
\begin{pmatrix}
  V_{11}(x) & V_{12}(x) & \ldots & V_{1d}(x) \\
  V_{21}(x) & V_{22}(x) & \ldots & V_{2d}(x) \\
  \vdots & \vdots & \vdots & \vdots \\
  V_{d1}(x) & V_{d2}(x) & \ldots & V_{dd}(x) \\
\end{pmatrix}.
$$
By applying the It\^o formula we find that
\begin{multline*}
\eps^{-2} V(X^\eps_{\varTheta_C\wedge
t})=\eps^{-2} V(x_0) +\int_0^{\varTheta_C\wedge
t}\eps^{-2}
\lef\nabla V(X^\eps_{s}), b(X^\eps_{s})\rig ds \\
 +\int_0^{\varTheta_C\wedge t}\eps^{-1}
\lef\nabla V(X^\eps_{s}), \sigma(X^\eps_{s})dB_s\rig
+\int_0^{\varTheta_C\wedge
t}\frac{1}{2}\trace\Big(\Psi(X^\eps_s)a(X^\eps
_s)\Big)ds.
\end{multline*}
We choose  $M^\eps_t=\int_0^{t}\eps^{-1} \lef\nabla V(X^\eps_{s}),
\sigma(X^\eps_{s})dB_s\rig$, which has  the variation process
$
\langle M^\eps\rangle_t=\int_0^{t}\eps^{-2} \lef\nabla V(X^\eps_{s}),
a(X^\eps_{s})\nabla V(X^\eps_s)\rig ds.
$
Clearly
\begin{multline*}
M^\eps_{\varTheta^\beta_C\wedge t}=\eps^{-2}
V(X^\eps_{\varTheta_C\wedge t})-\eps^{-2} V(x_0)
\\
-\int_0^{\varTheta_C\wedge t}\eps^{-2} \lef\nabla
V(X^\eps_{s}), b(X^\eps_{s})\rig ds
-\int_0^{\varTheta_C\wedge
t}\frac{1}{2}\trace\Big(\Psi(X^\eps_s)a(X^\eps
_s)\Big)ds.
\end{multline*}
Hence, by the definition of $\mathfrak{D}$,  one gets
\begin{multline}\label{hence}
M^\eps_{\varTheta_C\wedge T}-\frac{1}{2}\langle
M^\eps\rangle_{\varTheta_C\wedge T} =\eps^{-2}
V(X^\eps_{\varTheta_C\wedge T})-\eps^{-2} V(x_0)\\
-\int_0^{\varTheta_C\wedge T}
\frac{1}{2}\trace\Big(\Psi(X^\eps_s)a(X^\eps
_s)\Big)ds \quad -\int_0^{\varTheta_C\wedge T}
\eps^{-2}\mathfrak{D}V(X^\eps_s)ds.
\end{multline}
On the set $\{\varTheta_C\le T\}$, we have
$$
\eps^{-2} V(X^\eps_{\varTheta_C\wedge
T})-\eps^{-2} V(x_0)\ge
 \eps^{-2}\inf_{\|x\|\ge C}V(x)-\eps^{-2} V(x_0),
$$
and
$$
\Big|\int_0^{\varTheta_C\wedge T}
\frac{1}{2}\trace\big(\Psi(X^\eps_s)a(X^\eps
_s)\big)ds\Big|\le \frac{T}{2}\sup_{\|x\|\le
C}\big|\trace\big(\Psi(x)a(x)\big)\big|,
$$
and, by \eqref{a-2},
\begin{multline*}
-\int_0^{\varTheta\wedge T}
\eps^{-2}\mathfrak{D}V(X^\eps_s)ds
\\
\ge
-\Big|\int_0^{\varTheta_C\wedge T}
\eps^{-2}I_{\{\|X^\eps_s\|\le
L\}}\mathfrak{D}V(X^\eps_s)ds\Big|
\ge-\eps^2T\sup_{\|x\|\le L}|\mathfrak{D}V(x)|.
\end{multline*}
 These inequalities and \eqref{hence} imply
\begin{multline*}
M^\eps_{\varTheta_C}-\frac{1}{2}\langle
M^\eps\rangle_{\varTheta_C}\ge
\eps^{-2}\inf_{\|x\|\ge C}V(x)-\eps^{-2} V(x_0)
\\
-\frac{T}{2}\sup_{\|x\|\le C}\big|\trace\big(\Psi(x)a(x)\big)\big|
-\eps^{-2}T\sup_{\|x\|\le
L}\big|\mathfrak{D}V(x)\big|
\end{multline*}
on the set $\{\varTheta_C\le T\}$.
Hence, due to \eqref{PA-1} of Proposition \ref{pro-A.1}
\begin{multline*}
\eps^2\log\P\big(\varTheta_C\le T\big)\le
\\
-\inf_{\|x\|\ge C}V(x)+V(x_0)
+\frac{T\eps^2}{2}\sup_{\|x\|\le
C}\big|\trace\big(\Psi(x)a(x)\big)\big| +T\sup_{\|x\|\le
L}\big|\mathfrak{D}V(x)\big|
\\
\xrightarrow[\eps\to 0]{} -\inf_{\|x\|\ge C}V(x)+V(x_0)
+T\sup_{\|x\|\le
L}\big|\mathfrak{D}V(x)\big|
\end{multline*}
and it is left to recall that by \eqref{a-1} $\lim_{C\to\infty}\inf_{\|x\|\ge
C}V(x)=\infty$.
\end{proof}

\subsection{The proof of (\ref{3.1abs})}
We apply Lemma \ref{lem-2.1} to
$$
V(x)=\frac{c\|x\|^2}{1+\|x\|},
$$
with a positive parameter $c\le \frac{1}{K}$ for
$K$ from \eqref{tri} of Theorem \ref{theo-2.1}. The function
$V(x)$ is twice continuously differentiable and satisfies
\eqref{a-1}. It is left to show that $V(x)$ satisfies
\eqref{a-2} as well.

Direct computations give
$
\nabla V(x)=c\frac{(2+\|x\|)\|x\|}{(1+\|x\|)^2}\frac{x}{\|x\|}.
$
Denote
$$
r(x):=\frac{(2+\|x\|)\|x\|}{(1+\|x\|)^2}
$$
and notice that
$r(x)\le 1$. By assumption \eqref{dva} of Theorem \ref{theo-2.1}, one can
choose $L>0$  sufficiently large  so that $\lef
x,b(x)\rig<0 $ for any $\|x\|\ge L$. On the other hand, by assumption \eqref{tri}
of Theorem \ref{theo-2.1},
$ -1+\frac{c}{2}\frac{\lef x,a(x)x\rig}{\|x\| \ |\lef x,b(x)\rig|}\le-\frac{1}{2}
$ for $\|x\|\ge L$ and
\begin{align*}
\mathfrak{D}V(x)&= \Bigg(c\frac{r(x)}{\|x\|}\lef
x,b(x)\rig+\frac{c^2r^2(x)}{2}\frac{\lef x,
a(x)x\rig}{\|x\|^2}\Bigg)
\\
&=\Bigg( -c\frac{r(x)}{\|x\|}\big|\lef x,b(x)\rig\big|+
\frac{c^2r^2(x)}{2}\frac{\lef x,a(x)x\rig}{\|x\|^2}\Bigg)
\\
&=cr(x)\frac{|\lef x,b(x)\rig|}{\|x\|}\Bigg(-1+
\frac{c}{2}r(x)\frac{\lef x,a(x)x\rig}{\|x\| \ |\lef
x,b(x)\rig|}\Bigg)
\\
&\le cr(x)\frac{|\lef x,b(x)\rig|}{\|x\|}\Bigg(-1+
\frac{c}{2}\frac{\lef x,a(x)x\rig}{\|x\| \ |\lef
x,b(x)\rig|}\Bigg)
\\
&\le -\frac{1}{2}cr(x)\frac{|\lef x,b(x)\rig|}{\|x\|}
\end{align*}
 and \eqref{a-2} follows.
\qed

\subsection{The proof of (\ref{3.1abd})}

The obvious inclusion
\begin{multline*}
\Big\{\sup_{t\le
\triangle}\|X^{\eps}_{\theta+t}-X^{\eps}
_\theta\|>\eta\Big\}
\\
 \subseteq \Big\{\sup_{t\le
\triangle}\|X^{\eps}_{\theta+t}-X^{\eps}
_\theta\|>\eta, \ \varTheta_C=\infty\Big\}\bigcup
\Big\{\varTheta_C\le T\Big\}
\end{multline*}
reduces the proof to verifying
\begin{equation}\label{L32}
\varlimsup_{\triangle\to 0}\varlimsup_{\eps\to
0}\eps^2\log \sup_{\theta\le T} \P\Big(\sup_{t\le
\triangle}\|X^{\eps}_{\theta+t}-X^{\eps} _\theta\|>
\eta, \ \varTheta_C=\infty\Big)=-\infty
\end{equation}
for any fixed $C$. Indeed if \eqref{L32} holds, then
\begin{align*}
&\lim_{\triangle\to 0}\varlimsup_{\eps\to
0}\eps^2\log \sup_{\theta\le T} \P\Big(\sup_{t\le
\triangle}\|X^{\eps}_{\theta+t}-X^{\eps} _\theta\|>
\eta\Big)
\\
&\le \lim_{\triangle\to 0}\varlimsup_{\eps\to
0}\eps^2\log \sup_{\theta\le T} \P\Big(\sup_{t\le
\triangle}\|X^{\eps}_{\theta+t}-X^{\eps}
_\theta\|> \eta, \ \varTheta_C=\infty\Big) \\
& \hskip 3.0in \bigvee
\varlimsup_{C\to\infty}\varlimsup_{\eps\to 0}\eps^2
\log\P(\varTheta_C\le T\big)
\end{align*}
and, thus, \eqref{3.1abd} is implied by \eqref{L32} and \eqref{TC}.
So, it is left to check \eqref{L32} for any entry $x^\eps_t$
of $X^\eps_t$:
\begin{equation*}
\lim_{\triangle\to 0}\varlimsup_{\eps\to
0}\eps^2\log \sup_{\theta\le T} \P\Big(\sup_{t\le
\triangle}|x^{\eps}_{\theta+t}-x^{\eps} _\theta|>
\eta, \ \varTheta_C=\infty\Big)=-\infty.
\end{equation*}
A generic entry of $X^\eps_t$ satisfies
$$
x^\eps_t=x^\eps_0+\int_0^t\gamma^\eps_sds+\eps
m^\eps _t,
$$
where $\gamma^\eps_t$ is $\mathbf{F}^\eps$-adapted
continuous random process and $m_t$ is $\mathbf{F}^\eps$-continuous
martingale with $\langle m^\eps\rangle_t= \int_0^t\mu^\eps_sds$.
Since $b$ and $\sigma$ are locally Lipschitz continuous functions,
there is a constant
$l_C$, such that $|\gamma^\eps_{\varTheta_C\wedge t}|\le l_C$
and $\mu^\eps_{\varTheta_C\wedge t}\le l_C$.
Taking into account that
$$
\Big\{\sup_{t\le
\triangle}\Big|\int_\theta^{\theta+t}\gamma^\eps_sds\Big|\ge
\eta, \ \varTheta_C=\infty\Big\}\subseteq\big\{l_C\triangle\ge
\eta\big\}= \varnothing, \ \text{for $\triangle<\eta/l_C$},
$$
it is left to verify
\begin{equation*}
\lim_{\triangle\to 0}\varlimsup_{\eps\to
0}\eps^2\log \sup_{\theta\le T} \P\Big(\sup_{t\le
\triangle}|\eps m^{\eps}_{\theta+t}-\eps
m^{\eps} _\theta|> \eta, \ \varTheta_C=\infty\Big)=-\infty.
\end{equation*}
Due to the obvious inclusion
\begin{multline*}
\Big\{\sup_{t\le\triangle}|\eps
m^{\eps}_{\theta+t}-\eps m^{\eps}
_\theta|> \eta, \ \varTheta_C=\infty\Big\}= \\
\hskip 1.5in \Big\{\sup_{t\le \triangle}|\eps
m^{\eps}_{\varTheta_C\wedge(\theta+t)}- \eps
m^{\eps} _{\varTheta_C\wedge \theta}|> \eta, \
\varTheta_C=\infty\Big\}
\\
\hskip 1.5in \subseteq \Big\{\sup_{t\le \triangle}|\eps
m^{\eps}_{\varTheta_C\wedge(\theta+t)}- \eps
m^{\eps} _{\varTheta_C\wedge \theta}|> \eta\Big\},
\end{multline*}
we shall verify
\begin{equation*}
\varlimsup_{\triangle\to 0}\varlimsup_{\eps\to
0}\eps^2\log \sup_{\theta\le T} \P\Big(\sup_{t\le
\triangle}|\eps
m^{\eps}_{\varTheta_C\wedge(\theta+t)}- \eps
m^{\eps} _{\varTheta_C\wedge\theta}|> \eta\Big)=-\infty.
\end{equation*}
Notice  that $n^\eps_t:=\eps m^\eps_{\varTheta_C\wedge(\theta+t)}-
\eps m^\eps_{\varTheta_C\wedge\theta}$ is a continuous
martingale relative to $(\mathscr{F}^\eps_{\varTheta_C\wedge\theta+t})_{t\ge 0}$
(see e.g. Ch. 4, \S 7 in \cite{LSMar}) with
$
\langle n^\eps\rangle_t=\eps^2
\int_{\varTheta_C\wedge\theta}^{\varTheta_C\wedge(\theta+t)}
\mu^\eps_sds\le \eps^2l_Ct.
$
By the statement \eqref{PA-4} of Proposition \ref{pro-A.1},
$\P\big(\sup_{t\le\triangle}|n^\eps_t|\ge
\eta\big)\le 2e^{-\eta^2/ (2l_C\eps^2\triangle)},$ so that, $
\varlimsup_{\eps\to 0}\eps^2\log
\P\big(\sup_{t\le\triangle}|n^\eps_t| \ge
\eta\big)\le
-\frac{\eta^2}{2l_C\triangle}\xrightarrow[\triangle\to
0]{}-\infty. $ \qed

\section{\bf Local LDP upper bound}
\label{sec-5}

We start with the observation that \eqref{3.2abs} holds if for any
$T>0$
\begin{equation}\label{Lub}
\varlimsup_{\delta\to 0}\varlimsup_{\eps\to 0}\eps^2
\log\P\Big(\sup_{t\le T}
\|X^{\eps}_t-u_t\|\le\delta, \
\varTheta_C=\infty\Big)\le-J_T(u),
\end{equation}
since by the inclusion
\begin{gather*}
\Big\{\sup_{t\le T} \|X^{\eps}_t-u_t\|\le\delta\Big\}
\subseteq \Big\{\sup_{t\le T} \|X^{\eps}_t-u_t\|\le\delta,
\ \varTheta_C=\infty\Big\}\bigcup\Big\{ \varTheta_C\le T\Big\}
\end{gather*}
we have
\begin{align*}
&\varlimsup_{\delta\to 0}\varlimsup_{\eps\to
0}\eps^2 \log\P\Big(\sup_{t\le T}
\|X^{\eps}_t-u_t\|\le\delta\Big)
\\
&\le \varlimsup_{\delta\to 0}\varlimsup_{\eps\to
0}\eps^2 \log\P\Big(\sup_{t\le T}
\|X^{\eps}_t-u_t\|\le\delta, \ \varTheta_C=\infty\Big)
\bigvee \varlimsup_{\eps\to 0}\eps^2
\log\P(\varTheta_C\le T\big),
\end{align*}
and, by \eqref{TC}, the last term goes to $-\infty$ as $C\to\infty$.

We omit the standard proof for  $u_0\ne x_0$ or $du_t\not\ll dt$ (see, e.g.
\cite{DZ}).
The rest of the proof is split into several steps.

\subsection{$\b{u_0=x_0}$, $\b{du_t\ll dt}$,
$\b{\int_0^T\|\dot{u}_s\|^2ds< \infty}$ } Define the set
$$
\mathfrak{A}=\Big\{\sup_{t\le T}
\|X^{\eps}_t-u_t\|\le\delta, \ \varTheta_C=\infty\Big\}.
$$
With a continuously differentiable vector-valued function
$\lambda(s)$ of the size $d$, let us introduce a continuous local martingale
$
U_t=\int_0^{t}\lef\lambda(s),
\eps\sigma(X^\eps_s)dB_s\rig
$
and its martingale exponential
$
\mathfrak{z}_t=e^{U_t-0.5\langle U\rangle_t},
$
where
$$
\langle
U\rangle_t=\int_0^{t}\eps^2\lef\lambda(s),
a(X^\eps_s)\lambda(s)\rig ds.
$$
It is well known that $
\mathfrak{z}_t
$ is a continuous
positive local martingale, as well as a supermartingale. Consequently,
$\E\mathfrak{z}_T\le 1$ and, therefore,
\begin{equation}\label{cru}
1\ge \E I_{\{\mathfrak{A}\}}\mathfrak{z}_T.
\end{equation}
The required upper bound for $\P(\mathfrak{A})$ is obtained by
estimating $\mathfrak{z}_T$ from below on $\mathfrak{A}$.
Since
$ U_t=\int_0^{t}\lef\lambda(s), dX^{\eps}_s-
b(X^{\eps}_s)ds\rig$,
\begin{equation}\label{5.3bc}
\begin{aligned}
U_T-0.5\langle U\rangle_T&=\int_0^{T}\Big[\lef\lambda(s),
dX^{\eps}_s-
b(X^{\eps}_s)ds\rig
\\&\quad
-\frac{\eps^2}{2}\lef\lambda(s),
a(X^\eps_s) \lambda(s)\rig ds\Big]
\\
&=\int_0^{T}\Big[\lef\lambda(s), \dot{u}_s-b(u_s)\rig
-\frac{\eps^2}{2}\lef\lambda(s),
a(u_s)\lambda(s)\rig\Big]ds
\\
&\quad +\int_0^{T}\lef\lambda(s),
dX^{\eps}_s-\dot{u}_sds\rig
\\
&\quad +\int_0^{T}\lef\lambda(s), b(u_s)-b(X^{\eps}_s\rig
ds
\\
&\quad +\int_0^T\frac{\eps^2}{2}\lef\lambda(s),
[a(u_s)-a(X^{\eps}_s)] \lambda(s)\rig ds.
\end{aligned}
\end{equation}
We derive lower bounds on the set $\mathfrak{A}$ for each term in the right hand side of
\eqref{5.3bc}. Applying the It\^o
formula to $\lef \lambda(t), X^\eps_t-u_t\rig$, and taking
into account that $X^\eps_0=u_0$, we find that
\begin{gather*}
\lef \lambda(T), X^\eps_T-u_T\rig=\int_0^T\lef \lambda(s),
dX^\eps_s- \dot{u}_sds\rig+\int_0^T\lef \dot{\lambda}(s),
X^\eps_s- u_s\rig ds.
\end{gather*}
Therefore,
\begin{multline*}
\int_0^T\lef \lambda(s), dX^\eps_s- \dot{u}_sds\rig
\\
\ge -\Big| \lef \lambda(T), X^\eps_T-u_T\rig\Big|-
\Big|\int_0^T\lef \dot{\lambda}(s), X^\eps_s-u_s\rig
ds\Big|\ge -r_1\delta,
\end{multline*}
with $r_1:=r_1(\lambda,T,C)\ge 0$, independent of $\eps$.

Further, with  $r_i:=r_i(\lambda,T,C)\ge 0$, $i=2,3$, due to the local
Lipschitz continuity of $\sigma$ and $a$, we find that
\begin{align*}
& \int_0^{T}\lef\lambda(s), b(u_s)-b(X^{\eps}_s)\rig ds\ge
-r_2(\lambda,C,T)\delta
\\
& \int_0^T\frac{\eps^2}{2}\lef\lambda(s),
[a(u_s)-a(X^{\eps}_s)] \lambda(s)\rig ds\ge
-\eps^2r_3(\lambda,C,T)\delta.
\end{align*}

Hence with $r:=r_1+r_2+\eps^2 r_3$,
\begin{equation*}
\log\mathfrak{z}_T\ge \int_0^{T}\Big[\lef\lambda(s),
\dot{u}_s-b(u_s)\rig -\frac{\eps^2}{2}\lef\lambda(s),
a(u_s)\lambda(s)\rig\Big]ds -r(\lambda,T,C)\delta.
\end{equation*}
Set
$\nu(s)=\eps^2\lambda(s)$ and rewrite the above inequality
as:
\begin{multline*}
\log\mathfrak{z}_T\ge
\frac{1}{\eps^2}\int_0^{T}\Big[\lef\nu(s),
\dot{u}_s-b(u_s)\rig -\frac{1}{2}\lef\nu(s),
a(u_s)\nu(s)\rig\Big]ds
\\
-r\Big(\frac{\nu}{\eps^2},T,C\Big)\delta.
\end{multline*}
This lower bound, along with \eqref{cru}, provides the following
upper bound
\begin{multline*}
\eps^2\log\P \big(\mathfrak{A}\big)\le
-\int_0^{T}\Big[\lef\nu(s), \dot{u}_s-b(u_s)\rig
-\frac{1}{2}\lef\nu(s), a(u_s)\nu(s)\rig\Big]ds
\\
+\eps^2r\Big(\frac{\nu}{\eps^2},T,C\Big)\delta.
\end{multline*}
Clearly  $\varlimsup_{\eps\to
0}\eps^2r\big(\frac{\nu}{\eps^2},T, C\big)<\infty$
and, hence,
\begin{multline}\label{pochti}
\varlimsup_{\delta\to 0}\varlimsup_{\eps\to
0}\eps^2\log\P \big(\mathfrak{A}\big)
\\
\le
-\int_0^{T}\Big[\lef\nu(s), \dot{u}_s-b(u_s)\rig
-\frac{1}{2}\lef\nu(s), a(u_s)\nu(s)\rig\Big]ds.
\end{multline}
Since the left hand side of \eqref{pochti} is independent of
$\nu(s)$, \eqref{Lub} is derived by minimizing the right hand side
of \eqref{pochti} with respect  to $\nu(s)$. Two difficulties arise
on the way to direct minimization:
\begin{list}{-}{}
\item the matrix $a(u_s)$ may be singular
\item the entries of $\nu(s)$ should be continuously differentiable
functions.
\end{list}

\noindent Assume first $a(u_s)$ is a positive definite matrix, uniformly in $s$,
and write
\begin{multline*}
\lef\nu(s), \dot{u}_s-b(u_s)\rig -\frac{1}{2}\lef\nu(s),
a(u_s)\nu(s)\rig =\frac{1}{2}\|\dot{u}_s-b(u_s)\|^2_{a^{-1}(u_s)}
\\
-\frac{1}{2}\Big\|a^{1/2}(u_s)\big(\nu(s)-a^{-1}(u_s)[\dot{u}_s-b(u_s)]\big)\Big\|^2.
\end{multline*}
If the entries of $ a^{-1}(u_s)[\dot{u}_s-b(u_s)] $
are continuously differentiable functions, then, by taking
$\nu(s)\equiv -a^{-1}(u_s)[\dot{u}_s-b(u_s)]$ we find that
\begin{equation}\label{00}
\varlimsup_{\delta\to 0}\varlimsup_{\eps\to
0}\eps^2\log\P \big(\mathfrak{A}\big)\le
-\frac{1}{2}\int_0^{T}\|\dot{u}_s-b(u_s)\|^2_{a^{-1}(u_s)}ds.
\end{equation}
In the general case, due to $ \int_0^T\|\dot{u}_s\|^2ds<\infty, $ the
entries of  $a^{-1}(u_s)[\dot{u}_s-b(u_s)]$ are square integrable
with respect to the Lebesgue measure on $[0,T]$. Choose a
maximizing sequence $\nu_n(s)$, $n\ge 1$, of continuously
differentiable functions such that
$
\lim_{n\to\infty}\int_0^T\big\|\nu_n(s)-a^{-1}(u_s)[\dot{u}_s-b(u_s)]\big\|^2ds=0.
$
Since all the entries of $a(u_s)$ are uniformly
bounded on $[0,T]$
$$
\lim_{n\to\infty}\int_0^T\big\|a^{1/2}(u_s)\big(\nu_n(s)-a^{-1}(u_s)[\dot{u}_s-b(u_s)]
\big)\big\|^2ds=0
$$
and \eqref{00} holds too.

Now we drop the uniform nonsingularity assumption of
$a(u_s)$.
The upper bound in \eqref{00} remains valid with
$a(u_s)$ replaced by
$
a_\beta(u_s)\equiv a(u_s)+\beta \mathbf{I},
$
where $\beta$ is a positive number and $\mathbf{I}$ is $(d\times
d)$-unit matrix:
\begin{equation*}
\varlimsup_{\delta\to 0}\varlimsup_{\eps\to
0}\eps^2\log\P \big(\mathfrak{A}\big) \le
-\frac{1}{2}\int_0^T
\|\dot{u}_s-b(u_s)\|^2_{[a(u_s)+\beta\mathbf{I}]^{-1}}ds.
\end{equation*}

For any fixed $s$, the function $
\|\dot{u}_s-b(u_s)\|^2_{[a(u_s)+\beta\mathbf{I}]^{-1}} $ increases
with $\beta\downarrow 0$ and by Lemma \ref{lem-pi} possesses the limit
\begin{multline*}
\lim_{\beta\to 0}\|\dot{u}_s-b(u_s)\|^2_{[a(u_s)+
\beta\mathbf{I}]^{-1}}
\\
=
  \begin{cases}
 \|\dot{u}_s-b(u_s)\|^2_{a^\oplus(u_s)}, & \begin{array}{r}
a(u_s)a^\oplus(u_s)[\dot{u}_s-
 b(u_s)]
 \\=[\dot{u}_s-b(u_s)]
 \end{array} \\
    \infty , & \text{otherwise}.
  \end{cases}
\end{multline*}
Thus the required upper bound
\begin{multline*}
\varlimsup_{\delta\to 0}\varlimsup_{\eps\to
0}\eps^2\log\P \big(\mathfrak{A}\big)
\\
 \le
  \begin{cases}
    -\int_0^T
\frac{1}{2}\|\dot{u}_s-b(u_s)\|^2_{a^\oplus(u_s)}ds , &
\begin{array}{r} a(u_s)a^\oplus(u_s)[\dot{u}_s-
 b(u_s)]
 \\=[\dot{u}_s-b(u_s)], \ \text{a.s.}
\end{array}
\\
    \infty , & \text{otherwise}
  \end{cases}
\end{multline*}
follows by the monotone convergence theorem.

\subsection{$\b{u_0=x_0}$, $\b{du_t\ll dt}$,
$\b{\int_0^T\|\dot{u}_s\|^2ds=\infty}$ }

We emphasize that $du_t\ll dt$ on $[0,T]$ implies $\int_0^T\|\dot{u}_s\|ds<\infty$
and return to the upper bound from \eqref{pochti}. Since $b$ and $\sigma$ are locally Lipschitz, one can choose a constant $L$
(depending on $u(s)$), so that,
$|\lef \nu(s),b(u_s)
\rig|\le \|b(u_s)\|\|\nu(s)\| \le L\|\nu(s)\|$ and
$\lef\nu(s),a(u_s)\nu(s)\rig\le L\|\nu(s)\|^2$.
Then, \eqref{pochti} implies
\begin{equation*}
\varlimsup_{\delta\to 0}\varlimsup_{\eps\to
0}\eps^2\log\P \big(\mathfrak{A}\big)\le
-\int_0^{T}\Big[\lef\nu(s), \dot{u}_s\rig-L\|\nu(s)\|
-\frac{L}{2}\|\nu(s)\|^2\Big]ds.
\end{equation*}
Let $\nu_n(s)$ be a sequence of continuously differentiable functions, approximating the
bounded (for each fixed $p>0$) function $L^{-1}\dot{u}_sI_{\{\|\dot{u}_s\|\le p\}}$
in the sense that
$
\lim_{n\to\infty}\int_0^T\|\frac{1}{L}\dot{u}_sI_{\{\|\dot{u}_s\|\le
p\}}- \nu_n(s)\|^2ds=0.
$
Thus,
$$
\varlimsup_{\delta\to 0}\varlimsup_{\eps\to
0}\eps^2\log\P \big(\mathfrak{A}\big)
\le
-\frac{1}{2L}\underbrace{\int_0^{T}\big\|\dot{u}_s\big\|^2I_{\{\|\dot{u}_s\|\le
p\}}ds} _{\uparrow\infty \ \text{as $p\uparrow\infty$}}
+\underbrace{\int_0^T\big\|\dot{u}_s\big\|ds}_{<\infty}
\xrightarrow[p\to\infty]{}-\infty
$$
\qed

\section{\bf Local LDP lower bound. }
\label{sec-6}

If
$
\varlimsup_{\delta\to 0}\varlimsup_{\eps\to
0}\eps^2\log\P \big(\sup_{t\le \varTheta_C\wedge
T}|X^\eps_t-u_t|\le\delta\big)
\le -J_T(u)=-\infty,
$
then the corresponding lower bound in the local LDP is equal $-\infty$ too.
So in this section we examine the local LDP lower bound for $J_T(u)$
from \eqref{Lub}
when $J_T(u)<\infty$. The latter means that we
may restrict ourselves to analyzing test functions with the
properties:
\begin{equation}\label{6.2cd}
\begin{aligned}
& {\bf (i)}\quad u_0=x_0
\\
& {\bf (ii)}\quad du_t\ll dt
\\
& {\bf (iii)}\quad
a(u_t)a^\oplus(u_t)[\dot{u}_t-b(u_t)]=[\dot{u}_t-b(u_t)] \quad
\text{a.s.}
\\
& {\bf (iv)}\quad
\int_0^T\|\dot{u}_t-b(u_t)\|^2_{a^\oplus(u_t)}dt<\infty, \ \forall
\ T>0
\\
& {\bf (v)} \quad \int_0^T\|\dot{u}_t\|^2dt<\infty.
\end{aligned}
\end{equation}
Another helpful observation is that  \eqref{3.2abd} holds if for
any $C>0$
\begin{equation}\label{hr}
\varliminf_{\delta\to 0}\varliminf_{\eps\to
0}\eps^2\log\P \Big(\sup_{t\le \varTheta_C\wedge
T}\|X^\eps_t-u_t\|\le\delta\Big)
\ge -J_T(u)
\end{equation}
due to
$$
\Big\{\sup_{t\le \varTheta_C\wedge
T}\|X^\eps_t-u_t\|\le\delta\Big\}\subseteq \Big\{\sup_{t\le
T}\|X^\eps_t-u_t\|\le\delta\Big\}\bigcup\Big\{\varTheta_C\le
T\Big\}
$$
and \eqref{TC}.

\subsection{Nonsingular $\b{a(x)}$}
\label{sec-6.1} In this section, the matrix $a(x)$ is assumed to be
uniformly nonsingular in $x\in \Real$, in the sense that
$
a(x)\ge \beta\mathbf{I}
$
for a
positive number $\beta$.
Let
$
\lambda(s):=
\sigma^{-1}(X^\eps_s)\big[\dot{u}_s-b(X^{\eps}_s)\big]
$
and introduce a martingale $ U_t=\int_0^{\varTheta_C\wedge
t}\frac{1}{\eps} \lef\lambda(s),dB_s\rig $
and its martingale exponential
$\mathfrak{z}_t=e^{U_t-0.5\langle U\rangle_t}$,
$t\le T$, where $ \langle U\rangle_t =\int_0^{\varTheta_C\wedge
t}\frac{1}{\eps^2} \|\lambda(s)\|^2ds.$

By  ({\bf iv}) and ({\bf v}) of \eqref{6.2cd}, $\langle U\rangle_T\le \text{const.}$ and, hence,
$\E\mathfrak{z}_T=1$. We
use this fact in order to define a new probability measure
$\mathsf{Q}^\eps$ by
$d\mathsf{Q}^{\eps}=\mathfrak{z}_{T}d\P$. Since
$\mathfrak{z}_{T}$ is positive $\P$-a.s., $\P\ll
\mathsf{Q}^\eps$ as well and
$d\P=\mathfrak{z}^{-1}_{T}d\mathsf{Q}^{\eps}$.

We proceed with the proof of \eqref{hr} by applying
\begin{equation}\label{5.15ab}
\P(\widetilde{\mathfrak{A}})=\int_{\widetilde{\mathfrak{A}}}
\mathfrak{z}^{-1}_{T}d\mathsf{Q} ^{\eps}
\end{equation}
to the set
$
\widetilde{\mathfrak{A}}=\big\{\sup_{t\le\varTheta_C\wedge T}
\|X^{\eps}_t-u_t\|\le\delta\big\},
$
and estimating from below the right hand side in \eqref{5.15ab}.
In order to realize this program, it is convenient to have a semimartingale
description of
the process $X^{\eps}_{\varTheta_C\wedge t}$
under
$\mathsf{Q}^\eps$. Recall that the random process $B_{\varTheta_C\wedge
t}$ is a martingale under $\P$ with the variation process
$\langle B\rangle_{\varTheta_C\wedge t}\equiv(\varTheta_C\wedge t)\mathbf{I}$.
It is well known (see e.g. Theorem 2, Ch. 4,\S 5 in \cite{LSMar}) that
$B_{\varTheta_C\wedge t}$ is a continuous semimartingale under
$\mathsf{Q}^\eps$
with the
decomposition
$
B_{\varTheta_C\wedge t}=\widetilde{B}_{t}+A^B_{t},
$
where $\widetilde{B}_t$ is a martingale (under $\mathsf{Q}^{\eps}$)
with $\langle\widetilde{B}\rangle_t\equiv \langle B\rangle_{\varTheta_C\wedge t}$ and,
by the Girsanov theorem,
$$
A^B_{t}= \int_0^{\varTheta_C\wedge
t}\frac{1}{\eps}
\sigma^{-1}(X^\eps_s)[\dot{u}_s-b(X^{\eps}_s)]ds.
$$
In particular,
$$
X^{\eps}_{\varTheta_C\wedge t}=u_{\varTheta_C\wedge
t}+\eps \int_0^{\varTheta_C\wedge
t}\sigma(X^\eps_s)d\widetilde{B}_s, \ t\le T, \quad
\mathsf{Q}^\eps\text{-a.s.}
$$
As the next preparatory step we derive the semimartingale decomposition of $U_t$
under
$\mathsf{Q}^\eps$. As before, the continuous
martingale $U_t$ under $\P$ is
transformed to a semimartingale under
$\mathsf{Q}^\eps$:
$$
U_{t}=\widetilde{U}_{t}+A^U_{t}
$$
with continuous $\mathsf{Q}^\eps$-martingale
$\widetilde{U}_{t}$, having the variation process $\langle
\widetilde{U}\rangle_t\equiv\langle U\rangle_t$, $\P$- and
$\mathsf{Q}^\eps$-a.s., and a continuous drift $A^U_t\equiv \langle U\rangle_t$.

Thus, $ U_t=\widetilde{U}_t+\langle U\rangle_t, \ t\le T, \
\mathsf{Q}^\eps\text{-a.s.} $ and, thereby,
$
\mathfrak{z}^{-1}_T=e^{-\widetilde{U}_T- \frac{1}{2}\langle
U\rangle_T}.
$
Consequently, \eqref{5.15ab} is transformed to
\begin{align*}
\P(\widetilde{\mathfrak{A}})&= \int_{\widetilde{\mathfrak{A}}}
\exp\Big(-\widetilde{U}_{T}-\frac{1}{2}\langle U\rangle
_{T}\Big)d\mathsf{Q}^{\eps}
\\
&= \int_{\widetilde{\mathfrak{A}}}
\exp\Big(-\widetilde{U}_T-\frac{1}{2\eps^2}
\int_0^{\varTheta_C\wedge
T}\|\dot{u}_s-b(X^{\eps}_s)\|^2_{a^{-1}(X^\eps_s)}ds\Big)
d\mathsf{Q}^{\eps}.
\end{align*}
We are now in the position to derive a lower bound for the right
hand side.
Replacing $\widetilde{\mathfrak{A}}$ with a
smaller set $\widetilde{\mathfrak{A}}\cap\mathfrak{B}$, where $
\mathfrak{B}=\big\{\big|\eps^2 \widetilde{U}_{T}\big|\le \eta
\big\}$, write
$$
\P(\widetilde{\mathfrak{A}})\ge \int_{\widetilde{\mathfrak{A}}\cap\mathfrak{B}}
\exp\Big(-\frac{\eta}{\eps^2}-\frac{1}{2\eps^2}
\int_0^{\varTheta_C\wedge
T}\|\dot{u}_s-b(X^{\eps}_s)\|^2_{a^{-1}(X^\eps_s)}ds
\Big) d\mathsf{Q}^{\eps}.
$$
By the local Lipschitz continuity of
$b,\sigma$ and the uniform nonsingularity of $a(x)$,
$$
\Big|\|\dot{u}_s-b(X^{\eps}_s)\|^2_{a^{-1}(X^\eps_s)}
-\|\dot{u}_s-b(u_s)\|^2_{a^{-1}(u_s)}\Big|\le
l_C(\|\dot{u}_s\|+1)^2\delta, \ \delta\le 1,
$$
on the set
$\widetilde{\mathfrak{A}}\cap\mathfrak{B}$ for any  $s\le\varTheta_C\wedge T$.
Then,
\begin{align*}
\P(\widetilde{\mathfrak{A}}) & \ge \int_{\widetilde{\mathfrak{A}}
\cap\mathfrak{B}}
\exp\Big(-\frac{\eta}{\eps^2} -\frac{\delta
l_C}{\eps^2}\int_0^T(\|\dot{u}_s\|+1)^2ds
\\
& \hskip 4.0cm-\frac{1}{2\eps^2} \int_0^{\varTheta_C\wedge
T}\|\dot{u}_s-b(u_s)\|^2_{a^{-1}(u_s)}ds\Big)
d\mathsf{Q}^{\eps}
\\
& \ge \int_{\widetilde{\mathfrak{A}}\cap\mathfrak{B}}
\exp\Big(-\frac{\eta}{\eps^2} -\frac{\delta
l_C}{\eps^2}\int_0^T(\|\dot{u}_s\|+1)^2ds
\\
&  \hskip 4.0cm-\frac{1}{2\eps^2} \int_0^{T}
\|\dot{u}_s-b(u_s)\|^2_{a^{-1}(u_s)}ds\Big)
d\mathsf{Q}^{\eps}.
\end{align*}
Consequently,
$$
\varliminf_{\eps\to
0}\eps^2\log\P(\widetilde{\mathfrak{A}}) \ge -\eta-\delta
l_C\int_0^T(\|\dot{u}_s\|+1)^2ds -J_T(u)+ \varliminf_{\eps\to
0}\eps^2\log
\mathsf{Q}^{\eps}\big(\widetilde{\mathfrak{A}}\cap\mathfrak{B}\big).
$$
We prove now that $ \varliminf_{\eps\to
0}\eps^2\log
\mathsf{Q}^{\eps}\big(\widetilde{\mathfrak{A}}\cap\mathfrak{B}\big)=0 $
by showing
\begin{equation*}
\lim_{\eps\to 0}
\mathsf{Q}^{\eps}\big(\varOmega\setminus\widetilde{\mathfrak{A}}\big)=0
\quad \text{and}\quad \lim_{\eps\to 0}
\mathsf{Q}^{\eps}\big(\varOmega\setminus\mathfrak{B}\big)=0.
\end{equation*}
To this end, recall that
\begin{equation}\label{6.8z}
\begin{aligned}
\varOmega\setminus\widetilde{\mathfrak{A}}&=\Big\{\eps\sup_{t\le T}
\Big\|\int_0^{\varTheta_C\wedge t}\sigma(X^{\eps}_s)
d\widetilde{B}_s\Big\|>\delta\Big\}
\\
\varOmega\setminus\mathfrak{B}&=\Big\{\eps
\Big\|\int_0^{\varTheta_C\wedge
T}\sigma^{-1}(X^{\eps}_s)[\dot{u}_s-b(X^{\eps}_s)]
d\widetilde{B}_s\Big\|>\eta\Big\}.
\end{aligned}
\end{equation}
We verify \eqref{6.8z} componentwise. Let  $L^\eps_t$
denote any entry of $ \int_0^{\varTheta_C\wedge
t}\sigma(X^{\eps}_s) d\widetilde{B}_s $ or $
\int_0^{\varTheta_C\wedge
t}\sigma^{-1}(X^{\eps}_s)[\dot{u}_s-b(X^{\eps}_s)]
d\widetilde{B}_s. $ We show that
\begin{equation}\label{LLL}
\lim_{\eps\to
0}\mathsf{Q}^{\eps}\big(\eps\sup_{t\le T}
\big|L^\eps_t\big|>\delta\big)=0 \ \text{and} \
\lim_{\eps\to 0}\mathsf{Q}^{\eps}\big(\eps
\big|L^\eps_T\big|>\delta\big)=0.
\end{equation}
In both cases, $L^\eps_t$ is a continuous
$\mathsf{Q}^\eps$-martingale with $ \langle L^\eps\rangle_t
=\int_0^tg(s)ds \ $
and $ \int_\varOmega\int_0^Tg(s)dsd\mathsf{Q}^\eps<\infty.
$
Then \eqref{LLL} holds by Doob's  inequality:
\begin{equation*}
\lim_{\eps\to
0}\mathsf{Q}^{\eps}\big(\eps\sup_{t\le T}
\big|L^\eps_t\big|>\delta\big)\le
\frac{4\eps^2}{\delta^2}\int_\varOmega\int_0^{T}
g(s)dsd\mathsf{Q}^\eps \xrightarrow[\eps\to 0]{}0.
\end{equation*}

Now, for any fixed $\delta$ and $\eta$,
$$
\varliminf_{\eps\to
0}\eps^2\log\P(\widetilde{\mathfrak{A}}) \ge -\eta-\delta
l_C\int_0^T(\|\dot{u}_s\|+1)^2ds -J_T(u).
$$
The required lower bound
$$
\varliminf_{\delta\to 0}\varliminf_{\eps\to
0}\eps^2\log\P(\widetilde{\mathfrak{A}}) \ge -J_T(u)
$$
follows by taking
 $ \lim_{\eta\to 0}\lim_{\delta\to 0}$. \qed

\subsection{General $\b{a(x)}$}

This part of the proof requires perturbation arguments.
The idea is to use the already  obtained  local LDP
lower bound for the uniformly nonsingular $a(x)$.
Let $W_t$ be a standard $d$ dimensional Brownian motion, independent of  $B_t$,
defined on the same stochastic basis. Since $b$ and $\sigma$ are
assumed to be locally Lipschitz continuous, one can
introduce the perturbed diffusion process controlled by a free
parameter $\beta\in(0,1]$:
\begin{equation}\label{bbbb}
X^{\eps,\beta}_t=x_0+\int_0^tb(X^{\eps,\beta}_s)ds+\eps
\int_0^t[\sigma(X^{\eps,\beta}_s)dB_s+\sqrt{\beta}
dW_s].
\end{equation}
The process
$X^{\eps,\beta}_t$, defined in \eqref{bbbb},
solves the It\^o equation
$
X^{\eps,\beta}_t
=x_0+\int_0^tb(X^{\eps,\beta}_s)ds+\eps
\int_0^t[a(X^{\eps,\beta}_s)+\beta\mathbf{I}]^{1/2}dB^{\beta}_s.
$
with respect to a standard Brownian motion
$
B^{\beta}_t=\int_0^t[a(X^{\eps,\beta}_s)+\beta\mathbf{I}]^{-1/2}
[\sigma(X^{\eps,\beta}_s)dB_s+\sqrt{\beta} dW_s].
$
Then the family $\{(X^{\eps,\beta}_t)_{t\le
T}\}_{\eps\to 0}$ satisfies the local LDP lower bound.
Indeed, the matrix $ a_\beta(x)$ is uniformly
nonsingular, its entries are locally bounded and satisfy
the assumption \eqref{tri} of Theorem \ref{theo-2.1} since
$$
\frac{\lef x,a_\beta(x)x\rig}{\|x\| \ |\lef x,b(x)\rig|}
=\frac{\lef x,a(x)x\rig}{\|x\| \ |\lef x,b(x)\rig|}
+\beta\frac{\|x\|}{ |\lef x,b(x)\rig|}
$$
and $ \frac{\|x\|}{ |\lef
x,b(x)\rig|} $ converges to zero as $\|x\|\to\infty$ by
\eqref{dva}.
In
particular, with
$
\varTheta^\beta_C=\inf\{t:\|X^{\eps,\beta}_t\|\ge C\}
$
and $u_0=x_0, \ du_t\ll dt, \
\int_0^T\|\dot{u}_t\|^2dt<\infty$, we have
\begin{multline}\label{lbt}
\varliminf_{\delta\to 0}\varliminf_{\eps\to 0}
\eps^2\log\P\Big(\sup_{t\le \varTheta^\beta_C\wedge
T}\|X^{\eps,
\beta}_t-u_t\|\le\delta\Big)\ge \\
-\frac{1}{2}\int_0^T
\|\dot{u}_s-b(u_s)\|^2_{(a(u_s)+\beta\mathbf{I})^{-1}} ds.
\end{multline}
Further, we will use \eqref{lbt} to establish
\begin{equation}\label{1bto}
\varliminf_{\delta\to 0}\varliminf_{\eps\to 0}
\eps^2\log\P\Big(\sup_{t\le
T}\|X^\eps_t-u_t\|\le\delta\Big)\ge -\frac{1}{2}\int_0^T
\|\dot{u}_s-b(u_s)\|^2_{a^\oplus(u_s)} ds.
\end{equation}

To this end, we introduce the filtration $\mathbf{G}^\eps
=(\mathscr{G}^\eps_t)_{t\ge 0}$, with the general
conditions, generated by
$(X^\eps_t,X^{\eps,\beta}_t)_{t\ge 0}$ and notice
that both $\varTheta_C$ (see \eqref{ttt}) and $ \varTheta^\beta_C
$ are stopping times relative to $\mathbf{G}^\eps$. Hence,
\begin{equation}\label{tbc}
\tau^\beta_C=\varTheta_C\wedge \varTheta^\beta_C
\end{equation}
is a stopping time as well relative to
$\mathbf{G}^\eps$.
Obviously,
$$
\lim_{C\to\infty}\varlimsup_{\eps\to
0}\eps^2\log\P\big(\tau^\beta _C\le T\big)=-\infty.
$$
However, the proof of \eqref{1bto} requires a stronger
property:
\begin{equation}\label{str}
\lim_{C\to\infty}\varlimsup_{\eps\to
0}\eps^2\log\sup_{\beta\in(0,1]} \P\Big(
\tau^\beta_C\le T\Big)=-\infty.
\end{equation}

It is clear, that \eqref{str} is valid if it is valid with
$\tau^\beta_C$ replaced by $\varTheta^\beta_C$. The latter is
verified along the lines of Lemma \ref{lem-2.1} proof:
\begin{multline*}
\eps^2\log\sup_{\beta\in(0,1]}\P\big(\varTheta^\beta_C\le
T\big)\le -\inf_{\|x\|\ge C}V(x)+V(x_0)
\\
+\frac{T\eps^2}{2}\sup_{\beta\in(0,1]}\sup_{\|x\|\le
C}\big|\trace\big(\Psi(x)[a(x)+\beta\mathbf{I}]\big)\big|
+T\sup_{\beta\in(0,1]}\sup_{\|x\|\le L}\big|\mathfrak{D}_\beta V(x)\big|
\\
\xrightarrow[\eps\to 0]{} -\inf_{\|x\|\ge C}V(x)+V(x_0)
+T\sup_{\beta\in(0,1]}\sup_{\|x\|\le
L}\big|\mathfrak{D}_\beta V(x)\big|\xrightarrow[C\to\infty]{}-\infty,
\end{multline*}
where
$
\mathfrak{D}_\beta V(x)=\lef \nabla V(x), b(x)\rig+\frac{1}{2} \lef
\nabla V(x), a_\beta(x)\nabla V(x)\rig.
$
We are now in the position to prove \eqref{1bto}.
With $\delta\le \beta^{1/4}$, write
\begin{align*}
& \Big\{\sup_{t\le\tau^\beta_C\wedge
T}\|X^{\eps,\beta}_t-u_t\|\le\delta \Big\}
\\
&=\Big\{\sup_{t\le\tau^\beta_C\wedge
T}\|X^{\eps,\beta}_t-u_t\|\le\delta \Big\} \bigcap
\Big\{\sup_{t\le\tau^\beta_C\wedge T}
\|X^\eps_t-X^{\eps,\beta} _t\|\le\beta^{1/4}\Big\}
\\
&\quad\bigcup\Big\{\sup_{t\le\tau^\beta_C\wedge
T}\|X^{\eps,\beta}_t-u_t\| \le\delta\Big\} \bigcap
\Big\{\sup_{t\le\tau^\beta_C\wedge T}
\|X^\eps_t-X^{\eps,\beta}_t\|>\beta^{1/4}\Big\}
\\
&\subseteq \Big\{\sup_{t\le\tau^\beta_C\wedge
T}\|X^{\eps,\beta}_t-u_t\| \le\beta^{1/4}\Big\} \bigcap
\Big\{\sup_{t\le\tau^\beta_C\wedge T}
\|X^\eps_t-X^{\eps,\beta}_t\|\le\beta^{1/4}\Big\}
\\
&\quad\bigcup \Big\{\sup_{t\le\tau^\beta_C\wedge
T}\|X^\eps_t-X^{\eps,\beta} _t\|>\beta^{1/4}\Big\}
\\
&\subseteq \Big\{\sup_{t\le\tau^\beta_C\wedge
T}\|X^\eps_t-u_t\|\le 2\beta^{1/4} \Big\} \bigcup
\Big\{\sup_{t\le\tau^\beta_C\wedge
T}\|X^\eps_t-X^{\eps,\beta} _t\|>\beta^{1/4}\Big\}
\\
&\subseteq \Big\{\sup_{t\le T}\|X^\eps_t-u_t\|\le
2\beta^{1/4} \Big\}
\bigcup  \Big\{\sup_{t\le\tau^\beta_C\wedge
T}\|X^\eps_t-X^{\eps,\beta} _t\|>\beta^{1/4}\Big\} \\
& \hskip 8.7cm\bigcup
\Big\{\tau^\beta_C\le T\Big\}
\end{align*}
Hence,
\begin{multline*}
\P\Big(\sup_{t\le \tau^\beta_C\wedge
T}\|X^{\eps,\beta}_t-u_t\|\le\delta \Big) \le 3
\Bigg\{\P\Big(\sup_{t\le T}\|X^\eps_t-u_t\|\le
2\beta^{1/4} \Big)
\\
\bigvee \P\Big(\sup_{t\le\tau^\beta_C\wedge
T}\|X^\eps_t-X^{\eps,\beta} _t\|>\beta^{1/4}\Big)
\bigvee\P\big(\tau^\beta_C\le T\big)\Bigg\}.
\end{multline*}
Clearly,  $\varTheta^\beta_C$
can be replaced   by $\tau^\beta_C$, and so
\begin{multline}\label{6.16x}
- \frac{1}{2}\int_0^T
\|\dot{u}_s-b(u_s)\|^2_{(a(u_s)+\beta\mathbf{I})^{-1}}ds \le
\varliminf_{\eps\to 0}\eps^2
\log\P\Big(\sup_{t\le T} \|X^{\eps}_t-u_t\|\le
2\beta^{1/4}\Big)
\\
\bigvee \varlimsup_{\eps\to
0}\eps^2\log\P\Big( \sup_{t\le\tau^\beta_C\wedge
T}\|X^\eps_t-X^{\eps,\beta} _t\|>\beta^{1/4}\Big)
\\
\bigvee\varlimsup_{\eps\to
0}\eps^2\log\sup_{\beta\in(0,1]} \P\Big(
\tau^\beta_C\le T\Big).
\end{multline}
Recall the following facts:

1) by Lemma \ref{lem-pi} and \eqref{6.2cd},
$$
\lim_{\beta\to 0}\int_0^T
\|\dot{u}_s-b(u_s)\|^2_{(a(u_s)+\beta\mathbf{I})^{-1}}ds =\int_0^T
\|\dot{u}_s-b(u_s)\|^2_{a^\oplus(u_s)}ds;
$$

2) by Lemma \ref{lem-6.1a},
$$
 \lim_{\beta\to
0}\varlimsup_{\eps\to 0}\eps^2\log\P\Big( \sup_{t\le\tau^\beta_C\wedge
T}\|X^\eps_t-X^{\eps,\beta}
_t\|>\beta^{1/4}\Big)=-\infty;
$$

3) by \eqref{str}, $ \lim_{C\to\infty}\varlimsup_{\eps\to
0}\eps^2\log\sup_{\beta\in(0,1]} \P\Big(
\tau^\beta_C\le T\Big)=-\infty. $

\smallskip
\noindent
Hence, passing to the limit  $\beta\to 0$
and then $C\to \infty$ in  \eqref{6.16x} and taking into account 1) - 3),
one gets the required lower bound
$$
\varliminf_{\beta\to 0}\varliminf_{\eps\to 0}\eps^2
\log\P\Big(\sup_{t\le T} \|X^{\eps}_t-u_t\|\le
2\beta^{1/4}\Big)\ge -\frac{1}{2}\int_0^T
\|\dot{u}_s-b(u_s)\|^2_{a^\oplus(u_s)}ds.
$$
\qed

\appendix
\section{\bf Exponential estimates for martingales}
\label{sec-A} \label{App-A}
\begin{proposition}\label{pro-A.1}
{\rm (Lemma A.1 in \cite{GuiLip})} Let $M=(M_t)_{t\ge 0}$,
$M_t\in\Real$,  be a continuous local martingale with $M_0=0$ and
the predictable variation process $\langle M\rangle_t$ defined on
some stochastic basis with general conditions. Let $\tau$ be a
stopping time, $\alpha$ and $B$ positive constants and
$\mathfrak{A}$  some measurable set.

\begin{enumerate}
\renewcommand{\theenumi}{\alph{enumi}}
\item \label{PA-1} if $M_\tau-\frac{1}{2}\langle
M\rangle_\tau\ge\alpha$ on $\mathfrak{A}$, then $
\P(\mathfrak{A})\le e^{-\alpha}$;

\item \label{PA-2} if $ M_\tau\ge\alpha $ and $\langle
M\rangle_\tau \le B $ on $\mathfrak{A}$, then $
\P(\mathfrak{A})\le e^{-\frac{\alpha^2}{2B}}; $

\item \label{PA-3} $ \P(\sup_{t\le T}|M_t|\ge \alpha,
\langle M\rangle_T\le B)\le 2e^{-\frac{\alpha^2}{2B}}; $

\item \label{PA-4} $ \P(\sup_{t\le T}|M_t|\ge \alpha)\le
2e^{-\frac{\alpha^2}{2B}} \bigvee \P(\langle
M\rangle_T>B). $
\end{enumerate}

\end{proposition}

\section{\bf Pseudoinverse of nonnegative definite matrices}
Let $A^\oplus$ be the Moore-Penrose pseudoinverse matrix of $A$ (see \cite{Al}).
\begin{lemma}\label{lem-pi}
For $d\times d$ nonnegative definite matrix $A$ and $x\in\Real^d$,
$$
\lim_{\beta\to 0} \lef x, (A+\beta \mathbf{I})^{-1} x\rig=\begin{cases}
\|x\|^2_{A^\oplus}, & AA^\oplus x=x\\
\infty, &\text{otherwise}.
\end{cases}
$$
\end{lemma}
\begin{proof}
Let $S$ be an orthogonal matrix, $S^*S=\mathbf{I}$, such that
$
D:=S^*AS
$
is a diagonal matrix.
Then, due to $S^*(A+\beta\mathbf{I})S=D+\beta\mathbf{I}$, we have
$
S^*(A+\beta\mathbf{I})^{-1}S=(D+\beta\mathbf{I})^{-1}
$
and
$
S(D+\beta\mathbf{I})^{-1}S^*=(A+\beta\mathbf{I})^{-1}.
$
Write ($y:=S^*x$)
\begin{align*}
\lef x,(A+\beta \mathbf{I})^{-1} x\rig&=\lef x,S(D+\beta \mathbf{I})^{-1}S^* x\rig
=\lef S^*x,(D+\beta \mathbf{I})^{-1}S^*x\rig
\\
&=\lef y,(D+\beta \mathbf{I})^{-1}y\rig= \lef y,(D+\beta \mathbf{I})^{-1}DD^\oplus y
\rig
\\
&\quad+\lef y,(D+\beta \mathbf{I})^{-1}(\mathbf{I}-DD^\oplus)y\rig.
\end{align*}
Since
$
\lim_{\beta\to 0}(D+\beta\mathbf{I})^{-1}DD^\oplus=D^\oplus,
$
one gets
$$
\lim_{\beta\to 0}\lef y,(D+\beta\mathbf{I})^{-1}DD^\oplus y\rig=\|y\|^2_{D^\oplus}
=\|x\|^2_{A^\oplus}
$$
while
$
\lim_{\beta\to 0}\lef y,(D+\beta \mathbf{I})^{-1}(\mathbf{I}-DD^\oplus)y
\rig\ne \infty
$
only if $(\mathbf{I}-DD^\oplus)y=0$.
Since the latter condition is nothing but
$(\mathbf{I}-AA^\oplus)x=0$, the desired statement holds.
\end{proof}

\section{\bf Exponential negligibility of
$\b{X^{\eps,\beta}_{t}-X^\eps_{t}}$} \label{sec-B}

We start with the auxiliary result.
\begin{proposition}\label{pro-C.1}
Let $Y_t$ be a nonnegative continuous semimartingale defined on a stochastic basis
{\rm (}with general conditions{\rm ):}
\begin{multline}\label{R1}
Y_t=\int_0^{t}h_1(s)Y_sds+\eps\int_0^{t}h_2(s)Y_sdM'_s
\\
+\eps\sqrt{\beta}\int_0^{t}h_3(s)\sqrt{Y_s}dM''_s
+\eps^2\beta\int_0^th_4(s)ds,
\end{multline}
where $h_i(s), i=1,\ldots,4$, are bounded predictable processes and $M'_t$, $M''_t$
are continuous martingales, $d\langle M'\rangle_t=m'(t)dt$,
$d\langle M''\rangle_t=m''(t)dt$, $\langle M',M''\rangle_t\equiv 0$ with bounded
$m'(t)$ and $m''(t)$.
Assume that for any $T>0$ and $\beta>0$,
\begin{equation}\label{inpart}
\lim_{L\to\infty}\varlimsup_{\eps\to 0}\eps^2\log{P}\Big(\sup_{t\le T}
\sqrt{Y_t}>L\Big)=-\infty.
\end{equation}
Then, for any $T>0$,
$$
\lim_{\beta\to 0}\varlimsup_{\eps\to 0}\eps^2\log\P
\Big(\sup_{t\le T}\big|Y_t\big|>\beta^{1/4}\Big)=-\infty.
$$
\end{proposition}
\begin{proof}
Obviously $Y_t$ solves an integral equation
\begin{equation*}
Y_t=\mathcal{E}_t\int_0^t\mathcal{E}^{-1}_s\Big[\eps\sqrt{\beta}
h_3(s)\sqrt{Y_s}dM''_s+\eps^2\beta h_4(s)ds\Big],
\end{equation*}
where
$
\mathcal{E}_t=\exp\big(\int_0^t[h_1(s)-\eps^20.5h^2_2(s)]ds
+\int_0^t\eps h_2(s)dM'_s\big).
$
Let for definiteness $|h_i|\le r$, where $r$ is a constant.
Then, with $\eps\le 1$,
$$
\sup_{t\le T}|\log \mathcal{E}_t|\le T(r+0.5r^2)+\sup_{t\le T}\Big|\eps\int_0^t
h_2(s)dM'_s\Big|.
$$
Hence the random variable  $\sup_{t\le T}|\log \mathcal{E}_t|$ is bounded on the set
$$
\big\{\sup_{t\le T}\big|\eps\int_0^t
h_2(s)dM'_s\big|\le C\big\}.
$$
Moreover,
it is exponentially tight in the sense that
\begin{equation}\label{P0}
\lim_{C\to\infty}\varlimsup_{\eps\to 0}\eps^2\log\P\Big(
\sup_{t\le T}|\log \mathcal{E}_t|>C\Big)=-\infty.
\end{equation}
The latter is implied by
\begin{equation}\label{P1}
\lim_{C\to\infty}\varlimsup_{\eps\to 0}\eps^2\log\P\Big(
\sup_{t\le T}\big|\eps\int_0^t
h_2(s)dM'_s\big|>C\Big)=-\infty
\end{equation}
since the martingale $N_t=\eps\int_0^t
h_2(s)dM'_s$ has $\langle N\rangle_t=\eps^2\int_0^th^2(s)m'(s)ds$
and, with some positive number $r_1$, we have
$
\eps^2h^2(s)m'(s)\le \eps^2r_1.
$
Then, by taking into account that
$\P\big(\langle N\rangle_T>\eps^2r_1T\big)=0$ and
applying the statement \eqref{PA-4} of Proposition \ref{pro-A.1}, we obtain
$
\P\big(\sup_{t\le T}|N_t|>C\big)\le 2e^{-C^2/(2\eps^2 r_1T)}
$
providing \eqref{P1}.

Now we estimate  $\sup_{t\le T}|Y_t|$ on the set $\big\{\sup_{t\le T}\big|
\log\mathcal{E}_t\big|\le C\big\}$. Write
\begin{equation*}
\sup_{t\le T}|Y_t|\le e^CTr\eps^2\beta
+e^C\sup_{t\le T}\Big|\int_0^t\mathcal{E}^{-1}_s\eps\sqrt{\beta}
h_3(s)\sqrt{Y_s}dM''_s\Big|.
\end{equation*}
This upper bound and \eqref{inpart}, \eqref{P0} reduces the proof of Proposition
\ref{pro-C.1}
to:
\begin{multline*}
\lim_{\beta\to 0}\varlimsup_{\eps \to 0}\eps^2\log\P\Big(
\sup_{t\le T}\Big|\int_0^t\mathcal{E}^{-1}_s\eps\sqrt{\beta}
h_3(s)\sqrt{Y_s}dM''_s\Big|>\beta^{1/4},
\\
\ \sup_{t\le T}
\sqrt{Y_t}\le L,
\ \sup_{t\le T}|\log \mathcal{E}_t|\le C
\Big)=\infty
\end{multline*}
for any $C>0$ and $L>0$.
Introduce the martingale
$$
N''_t=\int_0^t\mathcal{E}^{-1}_s\eps\sqrt{\beta}
h_3(s)\sqrt{Y_s}dM''_s
\ \text{with} \
\langle N''\rangle_t=\int_0^t\mathcal{E}^{-2}_s\eps^2\beta
h^2_3(s)Y_sm''(s)ds
$$
and denote
$
\mathfrak{C}=\big\{\sup_{t\le T}
\sqrt{Y_t}\le L,
\ \sup_{t\le T}|\log \mathcal{E}_t|\le C\big\}.
$
With $r_2\ge h^2_3(s)Lm''(s)$, we find that
$$
\langle N''\rangle_T\le e^{2C}r_2T\eps^2\beta.
$$
Hence,
\begin{multline*}
\P\Big(\sup_{t\le T}|N''_t|>\beta^{1/4}, \mathfrak{C}\Big)=
\P\Big(\sup_{t\le T}|N''_t|>\beta^{1/4}, \
\langle N''\rangle_T\le e^{2C}r_2T
\eps^2\beta, \mathfrak{C}\Big)
\\
\le\P\Big(\sup_{t\le T}|N''_t|>\beta^{1/4}, \langle N''\rangle_T\le e^{2C}
r_2T\eps^2\beta\Big).
\end{multline*}
By \eqref{PA-3} of Proposition \ref{pro-A.1}
the latter term is bounded from above by
$$
2\exp\Big(\dfrac{\beta^{1/2}}{2e^{2C}r_2T\eps^2\beta}\Big).
$$
Then we obtain
$$
\varlimsup_{\eps\to \infty}\eps^2\log\P
\Big(\sup_{t\le T}|N''_t|>\beta^{1/4}, \mathfrak{C}\Big)\le
-\frac{1}{2e^{2C}r_2T\beta^{1/2}}\xrightarrow[\beta\to 0]{}-\infty.
$$
\end{proof}

We apply Proposition \ref{pro-C.1} in order to prove

\begin{lemma}\label{lem-6.1a}
For any $T>0$ and $C>0$,
$$
\lim_{\beta\to 0}\varlimsup_{\eps\to
0}\eps^2\log\P \Big(\sup_{t\le\tau^\beta_C\wedge
T}\big\|X^{\eps,\beta}_{t}-X^\eps_{t} \big\|
>\beta^{1/4}\Big)= -\infty.
$$
\end{lemma}
\begin{proof}

Recall that $X^\eps_t$ and $X^{\eps,\beta}_t$ solve
\eqref{1.1} and \eqref{bbbb} respectively and $\tau^\beta_C$ is
given in \eqref{tbc}.
Set $\triangle_t^{\eps,\beta}=X_{\tau^\beta_C\wedge
t}^{\eps,\beta} -X_{\tau^\beta_C\wedge t}^ {\eps}$.
By \eqref{1.1} and \eqref{bbbb},
\begin{multline*}
\triangle_t^{\eps,\beta}=\int_0^{\tau^\beta_C\wedge t}
\big(b(X_{\tau^\beta_C\wedge s}^{\eps,\beta})-
b(X_{\tau^\beta_C\wedge s}^{\eps})\big)ds+
\\
+\eps\int_0^{\tau^\beta_C\wedge
t}\big(\sigma(X_{\tau^\beta_C\wedge s} ^{\eps,\beta})-
\sigma(X_{\tau^\beta_C\wedge s}^{\eps})\big)dB_s
+\eps\sqrt{\beta}W_{\tau^\beta_C\wedge t}.
\end{multline*}
Due to the local Lipschitz continuity of $b$ and $\sigma$ and with
 $0/0=0$, the vector-valued and matrix-valued functions:
$$
f(s)=\frac{b\big(X_{\tau^\beta_C\wedge
s}^{\eps,\beta}\big)- b\big(X_{\tau^\beta_C\wedge
t}^{\eps}\big)}{\|\triangle_s^{\eps,\beta}\|}
\quad\text{and} \quad
g(s)=\frac{\sigma\big(X_{\tau^\beta_C\wedge s}
^{\eps,\beta}\big)- \sigma\big(X_{\tau^\beta_C\wedge
s}^{\eps}\big)}{\|\triangle_s^{\eps,\beta}\|}
$$
are well defined and their entries are bounded by a constant
depending on $C$.
Hence
\begin{equation*}
\triangle_t^{\eps,\beta}=\int_0^{\tau^\beta_C\wedge t}
\|\triangle_s^{\eps,\beta}\|f(s) ds
+\eps\int_0^{\tau^\beta_C\wedge
t}\|\triangle_s^{\eps,\beta}\| g(s)dB_s
+\eps\sqrt{\beta}W_{\tau^\beta_C\wedge t}.
\end{equation*}
Since $\|\triangle_t^{\eps,\beta}\|^2
=\lef \triangle_t^{\eps,\beta},\triangle_t^{\eps,\beta}\rig $, by the It\^o formula, we
find that
\begin{equation}\label{B.3}
\begin{aligned}
\|\triangle_t^{\eps,\beta}\|^2&=\int_0^t2\|\triangle_s^{\eps,\beta}\|
\lef \triangle_s^{\eps,\beta},f(s)\rig ds
\\
&+\eps\int_0^{\tau^\beta_C\wedge
t}2\|\triangle_s^{\eps,\beta}\| \lef
\triangle_s^{\eps,\beta}, g(s)dB_s\rig
\\
&+\eps\sqrt{\beta}\int_0^{\tau^\beta_C\wedge t}2
\lef\triangle_s^{\eps,\beta}, dW_s\rig
\\
&+\eps^2\int_0^{\tau^\beta_C\wedge
t}\|\triangle_s^{\eps,\beta}\|^2
\trace\big[g(s)g^*(s)\big]ds
\\
&+\eps^2\beta(\tau^\beta_C\wedge t)d.
\end{aligned}
\end{equation}
Now, by letting
$
\phi(s)=\frac{2\lef\triangle_s^{\eps,\beta},f(s)\rig}
{\|\triangle_s^{\eps,\beta}\|} $ and $
d\widehat{B}_s=\frac{2\lef\triangle_s^{\eps,\beta},g(s)dB_s\rig}
{\|\triangle_s^{\eps,\beta}\|}, $
we rewrite \eqref{B.3} as:
\begin{equation}
\label{C.7a}
\begin{aligned}
\|\triangle_t^{\eps,\beta}\|^2&=\int_0^{\tau^\beta_C\wedge
t} \|\triangle_s^{\eps,\beta}\|^2
\big\{\phi(s)
\\
&\quad +\eps^2\trace[g(s)g^*(s)]\big\}
ds
\\
&\quad
 +\eps\int_0^{\tau^\beta_C\wedge
t}\|\triangle_s^{\eps,\beta}\|^2 d\widehat{B}_s
\\
&\quad
+\eps\sqrt{\beta}\int_0^{\tau^\beta_C\wedge t}
\|\triangle^{\eps,\beta}_s\|\frac{2
\lef\triangle_s^{\eps,\beta}, dW_s\rig}{\|\triangle^{\eps,\beta}_s\|}\\
&
\quad
+\eps^2\beta(\tau^\beta_C\wedge t)d.
\end{aligned}
\end{equation}

With the notations
\begin{list}{-}{}
\item $Y_t=\|\triangle_t^{\eps,\beta}\|^2$
\item  $h_1(s)=I_{\{\tau^\beta_C\le s\}}\big\{\phi(s)+\eps^2\trace[g(s)g^*(s)]\big\}$
\item  $h_2(s)\equiv 1$
\item $h_4(s)=I_{\{\tau^\beta_C\le s\}}d$
\item  $M'_t=\widehat{B}_t$,
$m'(s)=\dfrac{4\lef\triangle_s^{\eps,\beta},g(s)g^*(s)
\triangle_s^{\eps,\beta}\rig}
{\|\triangle_s^{\eps,\beta}\|^2}$
\item $M''_t=\int_0^{\tau^\beta_C\wedge t}2
\dfrac{\lef\triangle_s^{\eps,\beta}, dW_s\rig}{\|\triangle_s^{\eps,
\beta}\|}$, \ $m''(s)\equiv 4$,
\end{list}
the equation \eqref{C.7a} is in the form of \eqref{R1}. Since $h_i(s)$,
$i=1,\ldots,4$ are bounded and
$
\sqrt{Y}_t\equiv\|X_{\tau^\beta_C\wedge
t}^{\eps,\beta} -X_{\tau^\beta_C\wedge t}^ {\eps}\|
\le \|X_{\tau^\beta_C\wedge
t}^{\eps,\beta}\|+\|X_{\tau^\beta_C\wedge t}^ {\eps}\|\le 2C,
$
i.e., \eqref{inpart} holds too,
the statement of the lemma follows from Proposition \ref{pro-C.1}.
\end{proof}

\end{document}